\documentclass[11pt]{article}
\usepackage[ansinew]{inputenc}
\usepackage{graphicx}
\usepackage{color}

\usepackage{enumerate,latexsym}
\usepackage{latexsym}
\usepackage{amsmath,amssymb}
\usepackage{graphicx}
\usepackage{times}

\newfont{\bb}{msbm10 at 11pt}
\newfont{\bbsmall}{msbm8 at 8pt}
\def\rth{\mathbb{R}^3}
\def\R{\mathbb{R}}
\def\B{\mathbb{B}}
\def\N{\mathbb{N}}

\def\D{\mathbb{D}}

\def\esf{\mathbb{S}}

\newcommand{\ben}{\begin{enumerate}}
\newcommand{\bit}{\begin{itemize}}
\newcommand{\een}{\end{enumerate}}
\newcommand{\eit}{\end{itemize}}

\newcommand{\ed}{\end{document}}

\newcommand{\wt}{\widetilde}
\newcommand{\HH}{\mbox{\bb H}}

\newcommand{\cd}{\mathcal{D}}
\newcommand{\intc}{\operatorname{Int}}

\def\a{{\alpha}}

\def\g{{\gamma}}
\def\G{{\Gamma}}

\def\de{{\delta}}
\def\be{{\beta}}
\def\ve{{\varepsilon}}
\def\F{{\cal F}}
\def\V{{\cal V}}

\def\D{{\cal D}}
\def\E{{\cal E}}
\def\a{{\alpha}}
\def\r{\mathbb{R}}
\def\b{\mathbb{B}}

\def\n{\mathbb{N}}

\newtheorem{theorem}{Theorem}[section]
\newtheorem{lemma}[theorem]{Lemma}

\newtheorem{remark}[theorem]{Remark}

\newtheorem{definition}[theorem]{Definition}
\newtheorem{conjecture}[theorem]{Conjecture}

\newenvironment{proof}{\smallskip\noindent{\it Proof.}\hskip \labelsep}
{\hfill\penalty10000\raisebox{-.09em}{$\Box$}\par\medskip}

\begin{document}

\begin{title}{Calabi-Yau domains in three manifolds}
\end{title}

\begin{author}{Francisco Martín\thanks{This research is partially
supported by MEC-FEDER Grant no. MTM2007 - 61775.}
\and William H. Meeks, III\thanks{This material is based upon
   work for the NSF under Award No. DMS -
   0703213. Any opinions, findings, and conclusions or recommendations
   expressed in this publication are those of the authors and do not
   necessarily reflect the views of the NSF.}}
   \end{author}
\maketitle


\begin{abstract}  We prove that for every smooth compact Riemannian three-manifold
$\overline{W}$ with nonempty boundary, there exists a smooth properly
embedded one-manifold
$\Delta \subset W=\intc(\overline{W})$, each of whose components
is a  simple closed curve and
such that the domain $\cd = W - \Delta$ does not admit any
properly immersed open surfaces with at least one annular end,
bounded mean curvature, compact boundary (possibly empty)
and a complete induced Riemannian metric.

\vspace{.2cm} \noindent   {\it 2000 Mathematics Subject
Classification.} Primary 53A10; Secondary 49Q05, 49Q10, 53C42.
\end{abstract}
\maketitle

\section{Introduction.}

A natural question in the global theory of minimal surfaces,
first raised by Calabi in 1965~\cite{calabi} and later revisited by
Yau~\cite{Yau, Yau-2}, asks whether or not there exists a complete
immersed minimal surface in a bounded domain $\cd$ in $\rth$. As is
customary, we will refer to this problem as the Calabi-Yau problem
for minimal surfaces. In 1996, Nadirashvili~\cite{na1} provided the
first example of a complete, bounded, immersed minimal surface in
$\rth$. However, Nadirashvili's techniques did not provide
properness of such a complete minimal immersion in any bounded
domain. Under certain restrictions on $\cd$ and the topology of an
open surface\footnote{We say that a surface is {\it open} if it is
connected, noncompact and without boundary.} $M$, Alarcón, Ferrer,
Martín, and Morales~\cite{density,  marmor1, marmor2, marmor3}
proved the existence of a complete, proper minimal immersion of $M$
in $\cd$.  Recently, Ferrer, Martin and Meeks~\cite{fmm} have given
a complete solution to the ``{\bf proper Calabi-Yau problem for
smooth bounded domains}" by demonstrating that for every smooth
bounded domain $\cd \subset \rth$ and for every open surface $M$,
there exists a complete proper minimal immersion $f\colon M\to \cd$;
furthermore, in~\cite{fmm}, they proved that such an immersion
$f\colon M\to \cd$ can be constructed so that for any two distinct
ends $E_1$, $E_2$ of $M$, the limit sets $L(E_1)$, $L(E_2)$ in
$\partial \cd$ are disjoint compact sets\footnote{ See
Definition~\ref{def:limit} for the definition of the limit set of an
end of a surface in a three-manifold.}.

In contrast to the above existence results, in this paper we prove
the existence of nonsmooth bounded domains $\cd$ in $\rth$, and more
generally, domains $\cd$ inside  any  Riemannian three-manifold, for
which some open surface
$M$ can not be properly immersed into $\cd$
as a complete surface with bounded mean curvature. In this case, we
will say that  $\cd$ is a {\bf Calabi-Yau domain} for $M$. The
result described in the next theorem generalizes the main theorem of
Martín, Meeks and Nadirashvili in~\cite{mmn} which demonstrates the
existence of nonsmooth bounded domains in $\rth$ which do not admit
any complete, properly immersed minimal surfaces with compact
boundary (possibly empty) and  at least one annular end.

\begin{theorem} \label{th}Let $\overline{W}$ be a smooth  compact
Riemannian three-manifold with nonempty boundary and let
$W=\intc(\overline{W})$. There exists a properly embedded
one-manifold $\Delta \subset W$ whose path components are smooth
simple closed curves, such that $\cd=W-\Delta$ is a Calabi-Yau
domain for any surface with compact boundary (possibly empty) and at
least one annular end. In particular, $\cd$ does not admit any
complete, noncompact,  properly immersed surfaces of finite
topology, compact boundary and constant mean curvature.
\end{theorem}

\section{Notation and the description of $\Delta$.}

Before proceeding with the proof of the main theorem, we fix some
notation.

\begin{enumerate}
\item $\b(R)=\{x\in \rth\mid |x|<R\} \;{\rm and}\; \b=\b(1)$.
\item $\overline{\b(R)}=\{x\in \rth \mid |x|\leq R\} \;{\rm and}\;
\overline{\b}=\overline{\b(1)}$.
\item $\esf^2(R)=\partial \b(R) \;{\rm and}\; \esf^2=\partial \b$.
\item For $p\in \rth$ and $\ve>0$, $\B(p,\ve)=\{x\in \rth\mid d(p,x)<\ve\}$
is the open ball of radius $\ve$ centered at $p$.
\item For $n\in \n$, $\b_n=\b(1-\frac{1}{2^n})\; {\rm and}\;
\esf_n^2=\partial \b_n$.
\item For any set $F\subset \rth$, the cone on $F$ is  $$C(F)=
\{x\in \rth\mid x=ta \;\;{\rm where}\; t\in(0,\infty) \;{\rm and}\; a\in F\}.$$
\item For any set $F\subset \rth$  and $\ve>0$, let $F(\ve)=\{x\in
\rth \mid d(x,F)\leq \ve\} $ be the closed $\ve$-neighborhood of
$F$, where $d$ is the distance function in $\rth$.
\end{enumerate}

In the proof of Theorem~\ref{th}, we will need the following
definition.

\begin{definition}\label{def:limit}
Let $f\colon M \to \cd$ be a proper immersion of  surface $M$ with
possibly nonempty boundary into an open domain $\cd$ contained in a
three-manifold $N$ with possibly nonempty boundary. The {\bf limit
set} of $M$ is $$L(M)=\bigcap_{\a\in I}(\overline{f(M) -
f(E_\a)}),$$ where $\{ E_{\a}\}_{\a \in I}$ is the collection of
compact subdomains of $M$ and the closure $\overline{f(M) -
f(E_\a)}$ is taken in $N$. The {\bf limit set $L(e)$ of an end $e$
of $M$} is defined to be  the intersection of the limit sets all
properly embedded subdomains of $M$ with compact boundary which
represent $e$. Notice that $L(M)$ and $L(e)$ are  closed sets of
$\partial \cd$, and so each of these limit sets  is compact when $N$
is compact.
\end{definition}

First we will prove Theorem~\ref{th} in the case $\overline{W}$ is
the smooth closed Riemannian ball  $\overline{\b}\subset \rth$.  In
this case, we will construct a properly embedded 1-manifold
$\Delta\subset \b$ with path components consisting of smooth simple
closed curves such that every proper immersion $f\colon
A=\esf^1\times [0,\infty)\to \b - \Delta$ of an annulus with a
complete induced metric has unbounded mean curvature; this result
will prove Theorem~\ref{th} in the special case
$\overline{W}=\overline{\b}$. The proof of the case of
Theorem~\ref{th} when $\overline{W} $ is a smooth Riemannian ball,
or more generally, an arbitrary compact smooth  Riemannian manifold
with nonempty boundary follows from straightforward modifications of
the proof of the $\b-\Delta$ case; these modifications are outlined
in the last paragraph of the proof.

The first step in the construction of $\Delta$ is to create a
CW-complex structure $\Lambda$ on the open ball $\B$. Consider the
boundary $\partial$ of the box $[-1,1]\times [-1,1]\times
[-1,1]\subset \rth$. The surface $\partial$ has a natural structure
of a simplicial complex ${\cal X}_1$ with faces ${\cal F}_1=\{F_1,F_2,...,
F_6\}$ contained in planes parallel to the coordinate planes, edges
$\E_1=\{E_1,E_2,\ldots, E_{12}\}$ and vertices ${\cal
V}_1=\{v_1,v_2,\ldots, v_8\}$. Let ${\cal X}_2$ denote the related
refined simplicial complex obtained from ${\cal X}_1$ by adding vertices
to the centers of each of the faces of ${\cal F}_1$ and to the
centers of each of the edges in $\E_1$, thereby obtaining  new
collections $\F_2, \E_2, \V_2$ of faces, edges, and vertices.  In
this subdivision each face of ${\cal F}_2$ corresponds to subsquare
in one the faces in ${\cal F}_1$ with four line segments, each of
length one. Note that $\F_2$ has $6\cdot4$ faces, $\E_2$ has
$2\cdot6\cdot 4$ edges and $\V_2$ has $6\cdot 4+2$ vertices.
Continuing inductively the refining of the complex ${\cal X}_2$,
produces at the $n$-th stage a simplicial complex ${\cal X}_n$ with
$6\cdot4^{n-1}$ square faces $\F_n$, $2\cdot 6\cdot 4^{n-1}$ edges
$\E_n$ and $6\cdot 4^{n-1}+2$ vertices $\V_n$.

We define the 1-skeleton $\G$ of $\Lambda$ as follows:
$$\G=\bigcup_{k=1}^{\infty}\left[C(\E_k)\cap
\esf_k^2\right]\cup\left[C(\V_k)\cap
\left(\overline{\b}_{k+1}-\b_k\right)\right],$$ where $C(\E_k)$
denotes the cone $C(\cup\E_k)$. Extend the proper 1-dimensional
$CW$-complex $\G\subset\B$ to a proper $2$-dimensional
$CW$-subcomplex $\Lambda'$ of $\Lambda$ as follows. The faces of
$\Lambda'$ are the spherical squares in $\esf_k^2-\G$, as $k$ varies
in $\n$, together with the set of flat rectangles
$C(\a)\cap(\B_{k+1}-\overline{\B}_k)$, where $\a$ is a 1-simplex in
$\G\cap\esf_k^2$, as $k$ varies  in $\n$ and $\a$ varies  in
$\G\cap\esf_k^2$, see Figure 1 below. Let ${\cal F}$ denote the set
of faces of $\Lambda$. Finally, $\B-\Lambda'$ contains an infinite
collection ${\cal G}=\{G_\a \}_{\a \in I}$ of components which have
the appearance of a cube  which is a radial product of a spherical
square in some $\esf^2_k-\G$ with a small interval of length
$2^{-(k+1)}$, together with the special component $\B(\frac{1}{2})$.
The set ${\cal G} $ is the set of 3-cells in
$\Lambda$, which completes the construction of the $CW$-complex
structure $\Lambda$ of $\b$.
\begin{figure}[htbp]
    \begin{center}
      \includegraphics[width=.5\textwidth]{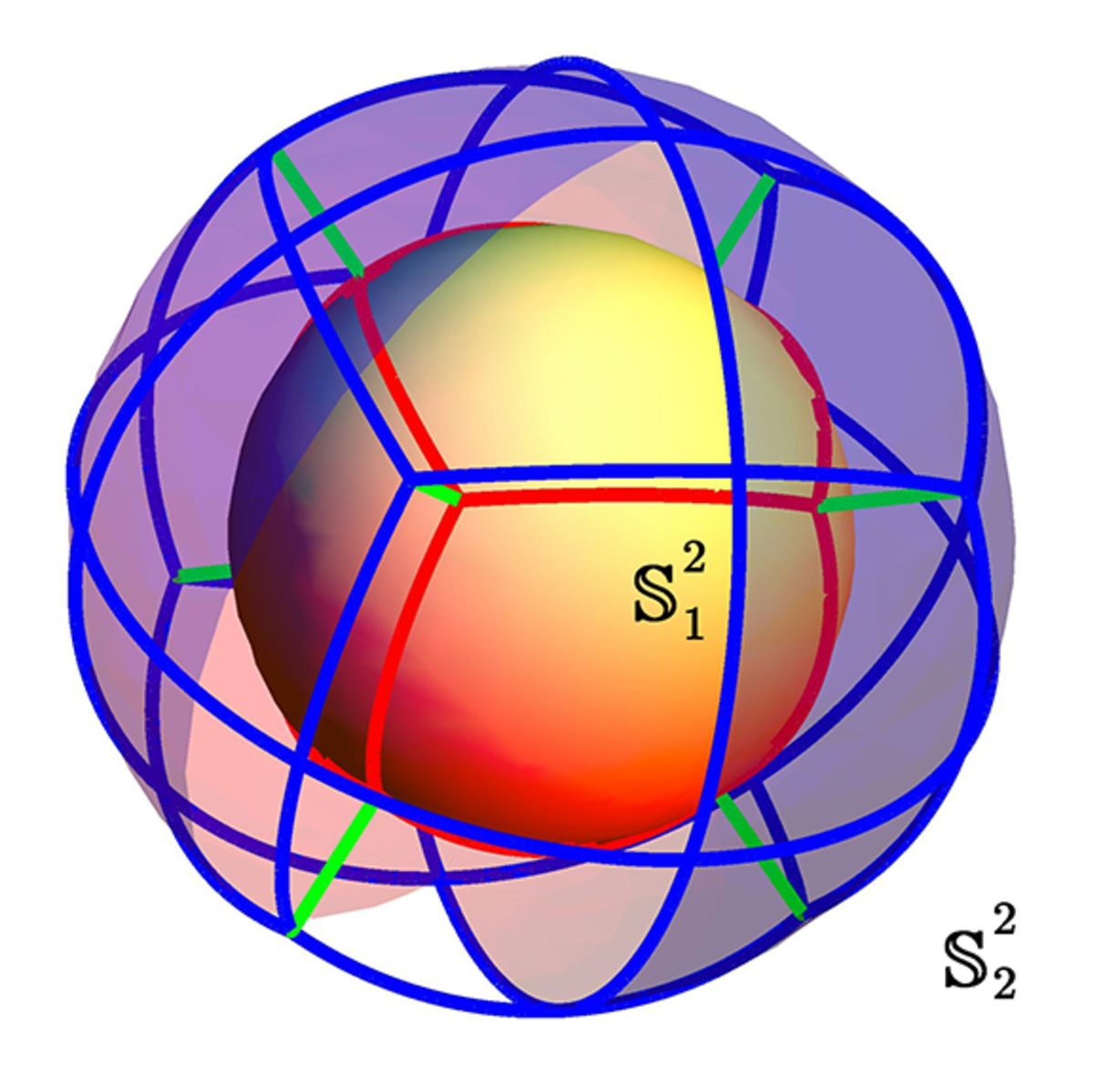}
   \end{center}
   \caption{The first two steps in the construction of $\Lambda$.}
\end{figure}

Define the related closed, piecewise smooth regular neighborhood
$\widehat{N}(\G)$ of $\G$:
$$\widehat{N}(\G)=\bigcup_{k=1}^{\infty}\left[\left(C(\E_k)\cap
\esf_k^2\right)\left(\frac{1}{2^k10}\right)\right]\cup
\left[\left(C(\V_k)\cap \left(
\overline{\b}_{k+1}-\b_k\right)\right)\left(\frac{1}{2^k
100}\right)\right].$$ Then let $N(\G)\subset \intc(\widehat{N}(\G))$
be a small smooth closed regular neighborhood of $\G$ in $\B$ such
that its boundary $\partial N(\G)$ intersects each face $F$ in
${\cal F}$ transversely in a simple closed curve $\be(F)$ that
bounds a disk $L(F)\subset F$; let ${\cal L} =\{L(F) \mid F \in
{\cal F}\}$.  For each open 1-simplex $\a \in \G$, let $P(\a)$ be
the plane perpendicular to $\a$ at the midpoint of $\a$. Let
$\wt{N}(\G)\subset \intc(\widehat{N}(\G))$ be another smooth closed
regular neighborhood of $\G$ with $N(\G)\subset \intc(\wt{N}(\G))$
and such that $\partial \wt N(\G) \cap P(\a)$ contains a simple
closed curve $\be(\a)$ close to $\a$ and which links $\a$. Let
$W(\a)\subset P(\a)$ denote the closed  disk with boundary curve
$\be(\a)$ and let ${\cal W}=\{W(\a) \mid \a \in \G\}$, see Figure 2.

\begin{figure}[htbp]
    \begin{center}
      \includegraphics[width=.75\textwidth]{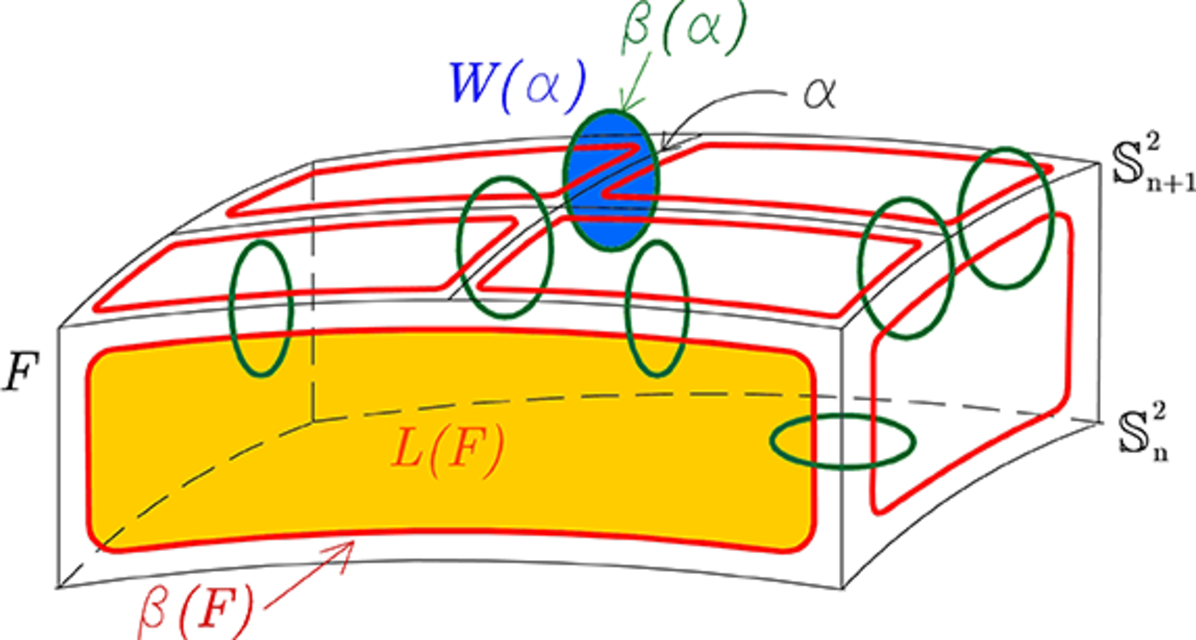}
   \end{center}
   \caption{The 1-dimensional simplicial complex $\G$,
   the 1-manifold $\Delta$ consisting of closed curves $\beta(F)$
   and  $\beta(F)$ and the disks $W(\a)$ and $L(F)$, where $F$ is a
   face in $\Lambda$ and $\a$ is a 1-simplex in $\Lambda$.}
\end{figure}

The set $\Delta$ is the collection
$\left[ \bigcup_{\a \in \G} \be(\a)\right] \cup \left[\bigcup_{F \in {\cal F}} \be(F)\right]$.
The domain described in Theorem~\ref{th} is ${\cal D}=\b-\Delta$.

We conclude this section with the following immediate consequence of
our constructions above.

\begin{lemma} \label{size}
Let $E$ be one of the following:\ben \item a 1-simplex, face or
3-cell in $\Lambda$;
\item a disk in either ${\cal W}$ or ${\cal L}$;
\item a component of $\wt{N}(\G) - \cup {\cal W}$. \een
If for some $\de \in (0,\frac{1}{4})$, $E \cap \left[\b -\b(1-\de)\right] \neq \text{\rm \O}$,
then $E$ is contained in a ambient ball $B_E$ of radius $4 \de$.
\end{lemma}

\section{$L(A)$ is a path connected subset of $\esf^2$ with more than one point.}
In this and the following sections, $f\colon A \to \D$ will denote  a
counterexample to Theorem~\ref{th} which, after a small smooth
perturbation, we will assume to be a fixed properly immersed annulus
diffeomorphic to $\esf^1 \times [0,1)$ satisfying:
\ben
\item The
supremum of the absolute mean curvature of $A$ is less than a fixed
constant  $H_0
>10$;
\item $f$ is transverse to the disks in ${\cal W}$ and to the
surface $\partial \wt{N}(\G)$;
\item $f$ is in general position with respect to $\Lambda$, i.e.,
$f$ is disjoint from the set of vertices ${\cal V}$ of $\Lambda$,
transverse to the closed faces of $\Lambda$ and so, it is also transverse
to $\esf^2_k$ for each $k\in \n$. \een
\begin{lemma}
\label{l3.1} If $f\colon\Sigma\to {\cal D}$ is a properly immersed
surface with compact boundary and $  {\bf e}$ is an end of $\Sigma$,
then the limit set $L({\bf e})$
of the end ${\bf e}$ is path connected.
\end{lemma}

\begin{proof}
This is a standard result, but for the sake of completeness, we
present its proof. Let $p,q\in L({\bf e})$ be distinct points. Let
$\D_1\subset \D_2\subset \ldots\subset
\D_n \subset\ldots$ be a smooth compact exhaustion of $\D$. After replacing
by subsequences, we may assume that there is a sequence of pairs of  points
${p}_n$, ${q}_n$ which
lie in the  component of $\Sigma -\intc(f^{-1}(\D_n))$ which
represents ${\bf e}$ and such that $\lim_{n \to \infty} f(p_n) = p$
and $\lim_{n \to \infty} f(q_n) = q$.

Let $\sigma_n\colon [0,1]\to \Sigma - \intc(f^{-1}(\D_n))$ be paths
with $\sigma_n(0)={p}_n$ and $\sigma_n(1)={q}_n$. Since the space
${\cal C}([0,1],\overline{\B})$ of continuous maps of $[0,1]$ into
$\overline{\B}$ is a compact metric space in the sup norm, a
subsequence of the paths  $f\circ\sigma_n$ converges to a continuous
map $f \circ \sigma$ of $[0,1]$ to  $\partial \D =\left[\esf^2\cup
\Delta\right] \subset \overline{\B}$ with $f\circ\sigma(0)=p$ and
$f\circ\sigma(1)=q$. Since $f\circ\sigma([0,1]) \subset L({\bf e})$
also holds, $L({\bf e})$ is path connected.
\end{proof}

\begin{lemma} \label{area1} If  $L(A) \cap \Delta \neq \mbox{\rm \O}$ or
if $L(A)$ consists of a single point in $\esf^2$, then $A$ has finite area.
\end{lemma}

\begin{proof}
By Theorems~3.1 and 3.1' in~\cite{hl1}, the bounded mean curvature
hypothesis and the properness hypothesis on $f$ imply that if $f$
composed with the inclusion map of $\D$ into $\rth$  is proper
outside of a point in $\esf^2$ or outside of a component of
$\Delta$, then the surface $A$ has finite area. Since $L(A)$ is
path connected and the path components of $\partial \D $ are
$\esf^2$ or a simple closed curve in $\Delta$, then the lemma
follows.
\end{proof}

\begin{lemma} \label{area2} If  $F\colon A \to \rth$ is a complete immersion
of \,$\esf^1 \times [0,\infty)$ with bounded mean curvature, then
$A$ has infinite area.
\end{lemma}
\begin{proof}
Suppose that $A$ has finite area and we will obtain a contradiction.
Since $A$ is a complete annulus of finite area, there exists a
sequence $\g_n$ of pairwise disjoint, piecewise smooth, closed embedded
geodesics with a single corner, which are topologically parallel to
$\partial A$ and whose lengths tend to $0$ as $n$ tends to infinity.
Assume that the index ordering of the geodesics $\g_n$ agrees with
the relative distances of these curves to $\partial A$. Replace $A$
by the subend $A(\g_1)$ with $\partial A(\g_1) =\g_1$. By the
Gauss-Bonnet formula applied to the subannulus $A(\g_1,\g_n)$ with
boundary $\g_1 \cup \g_n$, the total Gaussian curvature of
$A(\g_1,\g_n)$  is greater than $-4\pi$.  Since the Gaussian
curvature function $K_A$ of $A$ is pointwise bounded from above by
$H_0^2$, then the integral $\int_{A(\g_1)} K_A^+ \, dA$, where
$K_A^+(x) = \max \{K_A(x), 0\}$, is finite because $A$ has finite
area. Hence, after replacing $A$ by a subend of $A$, we may assume
that $\int_{A(\g_1)} K_A^+ \, dA <\pi$. So, we conclude that
$\int_{A(\g_1, \g_n)} K_A^- \, dA >-5\pi$, for all $n$, where
$K_A^-(x) = \min \{K_A(x), 0\}$.

On the other hand, since the area of $A$ does not grow at least
linearly with the distance from $\partial A$, the norm of the second
fundamental form of $A$ is unbounded on $A$. By standard rescaling arguments
(see for example~\cite{me29}), there exists a divergent
sequence $p_n\in A(\g_1)$ of blow-up points on the scale of the
second fundamental form with norm of the second fundamental form at
$p_n$ being $\lambda_n >n$, and intrinsic neighborhoods
$B_A(p_n,\frac{\lambda_n}{10})$ such that a subsequence of the
rescaled surfaces $\lambda_n\left[f(B_A(p_n,\frac{\lambda_n}{10}))-p_n\right]$
converges in the $C^2$-norm to  a  minimal disk $D$ in $\rth$ satisfying:
\ben \item The norm of the second
fundamental form of $D$ is at most 1 and equal to 1 at the origin.
\item $D$ is a graph over the projection to its tangent plane at the origin.
\item The total curvature of $D$ is  $-\ve$ for some $\ve >0$. Hence for
$n$ large, the integral of the function $K_A^-$ on
$B_A(p_n,\frac{\lambda_n}{10})$ is less than $-\frac{\ve}{2}$. \een
By property 3 above, we conclude that
$\lim_{n\to \infty}\int_{A(\g_1, \g_n)} K_A^- \, dA = -\infty$,
which contradicts our earlier observation that $\int_{A(\g_1, \g_n)}
K_A^- \, dA$ is bounded from below by $-5\pi$.
\end{proof}

The next lemma is an immediate consequence of Lemmas~\ref{area1}
and~\ref{area2}.

\begin{lemma} \label{l3.3} $L(A)$  is a path connected compact
subset of  $\esf^2$
containing  two distinct points $x$ and $y$. In particular, the
immersion $f$ can be
seen as a proper immersion in $\B$.
\end{lemma}

In the next sections, we will analyze how certain subdomains  of the
immersed annulus $f(A)$ intersects certain specific two-dimensional
subsets of $\cd$, for which we need the following definitions.

\begin{definition}{\rm Suppose $F\colon \Sigma \to \cd$ is a smooth
proper immersion of a surface with compact boundary which is
transverse to the disks in ${\cal W}$, to $\partial \wt{N}(\G)$ and
is in general position with respect to $\Lambda$. Suppose $\g$ is a
simple closed curve in $ \Sigma$. Then: \ben \item $\g $ is  an
\underline{$X_1$-type curve}, if $\g$ is a component of $F^{-1}(\cup
{\cal W})$.
\item $\g $ is an \underline{$X_2$-type
curve}, if $\g$ is a component of $F^{-1}(\partial \wt{N}(\G))$.
Note that in this case $\g \subset \left[\partial \wt N(\G)- \cup {\cal W}\right]$ and so
curves of $X_1$-type and $X_2$-type are disjoint.
\item $\g \subset  \Sigma$ is an \underline{$X_3$-type curve}, if $\g$ is a
component of $F^{-1}(\cup {\cal L})$. Notice that in this case
$\g$ is contained in a face of $\Lambda$. \een }
\end{definition}

\begin{definition}{\rm Given the fixed immersion $f\colon A \to \cd$, then:
\ben \item $X_1$ is the set of $X_1$-type curves  parallel to
$\partial A$ and  $X_2$ is the
set of $X_2$-type curves parallel to $\partial A$. \item   $X_3$ is the set of
$X_3$-type curves in  $A$ which are disjoint from $(\cup X_1) \cup
(\cup X_2).$\item By Lemma~\ref{lem:esf} below, the countable set $X$ can
be expressed as
$X= X_1 \cup X_2 \cup X_3 =\left\{\g_i\mid i\in \N\right\}$,
where the natural ordering of the simple closed,
pairwise-disjoint curves $\g_i$ in $A$ by their relative distances
from $\partial A$ agrees with the ordering of the index set $\N$.
\item $A_n$ denotes the compact subannulus in $A$ with $\partial
A_n=\partial A\cup \g_n$; note $A_1\subset A_2 \subset \ldots\subset
A_n \subset \ldots$ is a smooth compact exhaustion of $A$. \item For
$n,j\in\N$,  $A(n,j)$ denotes the compact subannulus of $A$ with
boundary curves $\g_n$ and $\g_{n+j}$.
\item $\displaystyle A(k)= \cup_{j=1}^\infty
A(k,j)$ is the end representative of $A$ with boundary $\g_k$.\een }
\end{definition}

\section{Placement properties of $\partial A(k,1)$ for k large.} \label{sec:plac}

\begin{lemma} \label{lem:esf}
For $k$ large, there exists at least one curve in $X$ in the region
$\B_{k+2} -\overline \B_{k-1}$.  In particular, the set $X$ is infinite.
\end{lemma}

\begin{proof} Assume that $f(\partial A)$ is contained in $\B_n$ and
we will prove that $\B_{k+2} - \overline \B_{k-1}$ contains an
element in $X$, whenever $k > n$.  Since $f \colon A \to \B$ is
proper and transverse to the spheres $\esf^2_i$ for every $i$, then
for $i\geq n$, $f^{-1}(\esf^2_i)$ contains a simple closed curve $\a_i$
which is a parallel to $\partial A$. If $f(\a_k) \cap (\cup
{\cal L}) \neq \O$, then either $\a_k \in X_3$ or  $\a_k$
intersects an element $\g$ of $X_1 \cup X_2$, where
$f(\g)$ is contained  in $\left[ \B_{k+1} - \overline \B_{k-1}\right] \subset \left[
\B_{k+2} - \overline \B_{k-1}\right]$. Similarly, if $f(\a_{k+1})\cap
(\cup{\cal L}) \neq \O$, then either $\a_{k+1} \in X_3$ or $\a_{k+1}$
intersects an element $\g$ of $X_1\cup X_2$, whose image $f(\g)$ must
be contained in $\left[
\B_{k+2} - \overline \B_{k}\right] \subset \left[ \B_{k+2} - \overline
\B_{k-1}\right]$. Hence, we may assume that $f(\a_k)$ and $f(\a_{k+1})$
are both disjoint from $\cup {\cal L}$ and so, $\left[f(\a_k \cup \a_{k+1})\right]
\subset \intc({\wt{N}}(\G))$.

Let ${D}_{\cal W}^k$ be the collection of disks in ${\cal W}$ which are
contained in $\B_{k+1} - \overline \B_k$ and let $\Sigma^k$ be the
compact domain which is closure of the component of $\wt{N}(\G) -
{(\cup D_{\cal W}^k})$ which contains $f(\a_k)$ in its interior. Let $A(\a_k,
\a_{k+1}) $ be the subannulus of $A$ with boundary $\a_k \cup
\a_{k+1}$.  Then $(f|_{A(\a_k, \a_{k+1})})^{-1}( \partial \Sigma^k) $
contains a simple closed curve $\g$ which is parallel to $\partial
A$ and which is an element of $X_1 \cup X_2 \subset X$. The existence of $\g$
completes the proof of the assertion.
\end{proof}

\begin{lemma} \label{lem:L2}
There exists a small $\eta_1>0$ such that for any
$\eta\in (0,\eta_1]$,  if $D \subset A$ is a compact
disk with $f(\partial D) \subset \B(z,\eta)$ for some $z \in \esf^2$ and $D$
contains a point $p$ such  that the distance $d(f(p),z)\geq1$, then:
\ben
\item  The disk $D$ contains
a $X_i$-type curve $\beta$, for $i=1,$ $2$ or $3$, and $f(\beta)$ lies
in $\B(z,1/2)- \overline{\B}(z,2 \eta)$.
\item The curve $\be$ can be chosen so that the disk
$D(\beta) \subset D$ bounded by $\beta$ contains $p$.
In particular, $f(D(\beta))$ contains a point of
distance at least $\frac12$ from its boundary and every
point in $D(\beta)$ has intrinsic distance at least $\eta$ from $\partial D$.
\een

\end{lemma}
\begin{proof}
Recall that for any face $F$ in ${\cal F}$, $C(F)$ denotes
the cone over $F$.  Clearly, for $\eta_1>0$ sufficiently small
and $\eta \in (0,\eta_1],$ there exist faces $F_1$, $F_2$, $F_3$
and $F_4$ in $\cal F$, such that:
$\B(z, 2 \eta) \subset \intc (C(F_1)) $, $ C(F_i) \subset \intc (C(F_{i+1}))$,
for $i=1,2,3$ and $C(F_4) \subset \B(z,1/2)$.

At this point we can follow the proof of Lemma~\ref{lem:esf}
where the annulus $D-\{p \}$ plays the role
of $A$ and the piecewise smooth disk $\partial
C(F_i)$ plays the role of $\esf^2_{k-2+i}$. Then we obtain an
$X_i$-type curve $\beta$ parallel to
$\partial D$ in $D-\{p\}$ and whose image $f(\beta)$ is in the
open region between $\partial (C(F_1))$
and $\partial (C(F_4))$, which is contained
$\B(z,1/2)- \overline{\B}(z,2\eta)$. This is the desired curve.
\end{proof}

Before stating the next assertion, we need some notation.

\begin{definition}
Given a curve $\g_k$ in $X$, we define $\chi_1(f(\g_k))$  to be  the
union of all closed $3$-cells in $\Lambda$ which intersect
$f(\gamma_k)$. Similarly, given $i \in \n$ we define
$\chi_{i+1}(f(\g_k))$ as the union of all closed $3$-cells in
$\Lambda$ which intersect $\chi_i\left(f(\gamma_k)\right)$.
\end{definition}
In what follows, we shall use the observation  that for
$i=1$ and $2$,  the set $\chi_i(f(\g_k))$ is a piecewise
smooth compact ball, whose boundary sphere is a union of
faces in $\cal F$ and it is in general position with respect to
the immersion $f$.

\begin{lemma} \label{lem:chi} For $k$ large, we have
$f(\gamma_{k+1}) \subset \chi_3(f(\gamma_k))$ or
$f(\gamma_{k}) \subset \chi_3(f(\gamma_{k+1}))$.
Furthermore, given $\eta>0$, there exists an integer
$k(\eta)$ such that for any  $k\geq k(\eta)$ one has:

\begin{enumerate}
\item $f(A(k)) \subset \left[\B -\overline\B(1-\eta)\right]$ and
each $X_i$-type curve $\g$, $i=1, 2$ or $3$, in $A(k,1)$ is
contained in a ball $\B(y(\g),\eta)$ for a suitable point $y(\g)\in \esf^2$.
\item There is a point $z(k)\in \esf^2$ such that
$f(\g_k \cup \g_{k+1} )\subset \B (z(k), \eta)$.
\item Every simple closed curve
$\g\subset \left[A(k,1)-f^{-1}(\B(z(k),\eta))\right]$ bounds a disk in $A(k,1)$.
\end{enumerate}

\end{lemma}
\begin{proof} In order to prove the first statement of the lemma,
we distinguish four cases, depending on the position
of $f(\gamma_k)$.  We will use the fact that by
Lemma~\ref{size}, for $k \to \infty$, the curve
$f(\g_k)$ becomes arbitrarily close to a point $z(k)\in \esf^2$.

\noindent {\bf Case A:} $f(\gamma_k) \subset D \in {\cal W}$.

In this case $f(A(k,1))$ enters a component $C$ of $\wt N - \cup
{\cal W}$ near $f(\gamma_k)$. Consider the compact component $Z$ of
$\left( f|_{A(k)}\right)^{-1}(\overline{C}) \subset A(k)$
with boundary component $\g_k$
and let
$\a_k$ be the boundary curve of $Z-\g_k$ which is parallel to
$\partial A(k)=\g_k$. By the  definition of $X_1$ and $X_2$,
$\a_k=\gamma_{k+j} \in \left[X_1 \cup X_2 \right]
\subset X$, for some $j \geq1$. By definition of $X$,
$\gamma_{k+1} \subset A(k,j)$ and intersects the domain $Z$. If
$\gamma_{k+1} \subset Z$, then clearly $f(\gamma_{k+1}) \subset
\overline{C} \subset \chi_2(f(\gamma_k))$ and we are done.
Otherwise, $f(\gamma_{k+1})$ must not be contained in $\wt
N(\Gamma)$. This means that $\g_{k+1}$ belongs to $X_3$ and so it is
contained in a face $F$ of $\Lambda$, which intersects $C$. Hence,
$F \subset \chi_2(f(\g_k))$ which implies $f(\g_{k+1}) \subset
\chi_2(f(\g_k))$.

\noindent {\bf Case B:} $f(\gamma_k) \subset \left[ \partial
\wt N(\G)- \cup {\cal W} \right]$ and the annulus
$f(A(k,1))$ enters $\widetilde N(\G)$ near $f(\g_k)$.

In this case, the arguments in Case A apply to show that
$f(\g_{k+1}) \subset \chi_2(f (\g_k))$.

\noindent {\bf Case C: } $f(\gamma_k) \subset
\left[ \partial \wt N(\G)- \cup {\cal W} \right]$ and the annulus
$f(A(k,1))$ enters $\B - \widetilde N(\G)$ near $f(\g_k)$.

%

First, note that if $f(\g_{k+1})$ intersects $\chi_2(f(\gamma_k))$,
then $f(\g_{k+1}) \subset \chi_3(f(\gamma_k))$. Thus, we may assume
that $f(\g_{k+1})$ lies outside the compact piecewise smooth ball
$\chi_2(f(\g_k))$. Consider the compact component $Z$ of
$(f|_{A(k,1)})^{-1}(\chi_2(f(\g_k)))$ containing $\g_k$ in its boundary.
Let $\a_k\neq \g_k$
be the boundary curve of $Z$ which is parallel in $A(k)$ to
$\g_k$; recall that $A(k)$ is the end of $A$ with boundary $\g_k$.
If $f(\a_k)$ intersects $\cup {\cal L}$, then $f(\a_k)$ is
contained in a disk $D \in {\cal L}$; in this case, since $\a_k$ lies between
$\g_k$ and $\g_{k+1}$ and it is parallel to $\partial A(k)$, then
$\alpha_k \in X_3$, which is contrary to the definition of
$\g_{k+1}$. Thus, $f(\a_k)\subset \partial(\chi_2(f(\g_k)))$ and  is disjoint
from $\cup{\cal L}$, and so $f(\a_k) \subset \intc (\wt N(\G))$. Let
$A(\g_k,\a_k) \subset A(k,1)$ be the subannulus with boundary curves
$\g_k \cup \a_k$. As $f(A(\g_k,\a_k))$ enters $\B - \wt N(\G)$ nears
$f(\g_k) $ and $f(\a_k) \subset \intc (\wt N(\G))$, then ours
previous separation arguments imply that there exists a curve
$\beta \subset \left(f \left|_{A(\g_k,\a_k)}
\right. \right)^{-1} (\partial \wt N(\G)- \cup {\cal W}) $ which is
parallel to $\g_k$. Since $\beta \in X_2$ and $\beta \neq \g_{k+1}$,
we arrive at a contradiction. This contradiction proves Case C.

\noindent {\bf Case D:} $f(\gamma_k) \subset D \in {\cal L}$.

If $f(\g_{k+1}) \subset \partial \wt N(\G)$ or $f(\g_{k+1}) \subset
\widehat D \in {\cal W}$, then the arguments in our previously
considered cases imply that $f(\g_k) \subset \chi_3(f(\g_{k+1}))$.
Hence, we may assume that $f(\g_{k+1}) \subset D' \in {\cal L}$ as
well.

If $\chi_1 (f(\g_k)) \cap \chi_1 (f(\g_{k+1})) \neq \O$, then
$f(\g_{k+1}) \subset \chi_2 (f(\g_k))$. Hence, we can assume that
$\chi_1 (f(\g_k)) \cap \chi_1 (f(\g_{k+1})) = \O$. Recall that
$f|_{A(k,1)}$ is in general position with respect to $\partial\left(
\chi_1 (f(\g_k)) \right) $ and $\partial \left(\chi_1 (f(\g_{k+1}))
\right)$. Let  $Z_i$ be the component of
$(f|_{A(k,1)})^{-1} (\chi_1(f(\g_i)))$ with boundary
component $\g_i$ and let $\a_i \neq \g_i$ be the boundary component of
$Z_i$ which is parallel to $\g_i$,
for $i=k, k+1$, respectively. Since $\a_k$ and $\a_{k+1}$ lie in
$\intc(A(k,1))$, then by definition of $X$, both $f(\a_k)$ and
$f(\a_{k+1})$ are disjoint from $\cup{\cal L}$. Moreover, as $f(\a_i)
\subset \chi_1 (f(\g_i))$, for $i=k,$ $k+1$, then $f(\a_k \cup
\a_{k+1}) \subset \intc ( \wt N(\G))$. Let $A(\a_k,\a_{k+1})$ be
the subannulus of $A(k,1)$ with boundary $\a_k \cup \a_{k+1}$.

Consider the collection of disks $ D^k_{\cal W}$ in $\cal W$ which are
contained in the interior of $\chi_2(f(\g_k)) - \chi_1 (f(\g_k)). $
Then $\wt N(\G) - \cup { D}^k_{\cal W}$ contains a connected  domain
whose closure $\Sigma^k$ in $\B$ satisfies
$f(\a_k) \subset \intc (\Sigma^k)$ and $f(\a_{k+1})
\subset \B-\Sigma^k$. Our previous separation arguments imply that
there is a simple closed curve $\beta$ in
$(f|_{A(\a_k,\a_{k+1})})^{-1} (\partial \Sigma^k)$ which is parallel
to $\g_k$. But $\beta \subset \intc(A(k,1)$ and
$\beta \in X_1 \cup X_2$, which is a contradiction. This
contradiction completes the proof of the first statement of the
lemma.

Item 1 in the lemma is a straightforward consequence of the fact that,
as $ \to \infty$, then $f(A(k))$ uniformly converges to $\esf^2$. Moreover,
given a $X_i$-type curve $\g \subset A(k)$, $i=1,2,3$,
the Euclidean diameter of $f(\g)$ goes to zero (as $k \to \infty$)
and is arbitrarily close to
a point $y(\g)$ in $\esf^2$.
Item  2 in the lemma follows from the observation that as
$k\to \infty$, the sets $\chi_3$ $(f(\g_k))$ are arbitrarily close
to $f(\g_k)$, which in turn, lie arbitrarily close to points
$z(k)\in \esf^2$.  These observations imply that there exists an integer
$j(\eta)$ such that for $k\geq j(\eta)$, items 1 and 2 in Lemma~\ref{lem:chi} hold.

In order to obtain item 3, we define $k(\eta)=j(\frac{\eta}{900})$. By definition of
$j(\frac{\eta}{900})$, for $k\geq k(\eta)$,
$f(\g_k\cup \g_{k+1})\subset \B(z(k),\frac{\eta}{900})$ and $f(A(k))\subset \left[\B-\overline{\B}(1-\frac{\eta}{900})\right]$. It remains to
check that each simple closed curve $\beta$ in a component $K$ of
$(f_{A(k,1)})^{-1}(\B-\B(z(k),\eta))$ bounds a disk in $A(k,1)$;
note that $\overline{K}\subset \intc(A(k,1))$. Observe that
$\frac{\eta}{900}$ is sufficiently small so that there exist
faces $F_1,F_2,F_3$ and $F_4$ in ${\cal F}$, such that:
$\B(z(k),\frac{\eta}{900})\subset \intc(C(F_1)), C(F_i)\subset \intc(F_{i+1}),$
for $i=1,2,3$ and $C(F_4)\subset \B(z(k),\eta)$.

If $\beta \subset K$ does not bound a disk in $A(k,1)$, then it is parallel to
$\g_k$ in $A(k,1)$. Let
$A(\g_k, \beta)$ denote the subannulus of $A(k,1)$ with boundary $\g_k\cup\beta$.
Then the arguments in the proof of Lemma~\ref{lem:L2} imply that there exists a
simple closed curve $\g'\subset \intc(A(\g_k, \beta))$ which is parallel to $\g_k$,
$f(\g')\subset \B(z(k), \eta)$ and $\g'$ is an $X_i$-type curve, for $i=1, 2$ or
$3$. In particular, $\g'\in X$ which is impossible. Thus, every simple closed curve
in $A(k,1)$ whose image under $f$ lies outside  of $\B(z(k), \eta)$ bounds a disk
in $A(k,1)$. This completes the proof of the lemma.
\end{proof}

The next lemma directly follows from the mean curvature comparison
principle.

\begin{lemma} \label{lem:maximum}
Suppose $\Sigma \subset A$ is a compact domain such that:
$f(\partial \Sigma)$ is contained in $\B(z, \eta)$, where $z \in \esf^2$
and $\eta<\frac{1}{H_0}$. Then either $f(\Sigma)\subset \B(z, \eta)$ or
$f(\Sigma)$ contains a point outside of $\B(z, \frac{1}{H_0})$.
\end{lemma}

\section{Proof of the Theorem~\ref{th}.}
By Lemma~\ref{l3.3}, $L(A) \subset\esf^2$ contains at least two
distinct points $x$ and $y$. We next prove that the limit set of $f$
is the entire sphere $\esf^2$.

\begin{lemma} \label{l4.1}
$L(A)=\esf^2$.
\end{lemma}

\begin{proof}
By Lemma~\ref{l3.3}, there are distinct points $x,y\in L(A)\subset
\esf^2$. Arguing by contradiction, suppose that there exists a point
$p\in \esf^2-L(A)$. The definition of limit point and the fact that
$f\colon A\to \D$ is proper with $L(A)\subset \esf^2$ imply  there
exists an $\ve>0$ such that $\B(p,10\ve)\cap f(A)=\text{\O}$. By
properness of $f$ in $\B$, then for $n$ large, we have
$f(\overline{A-A_n})\subset \left[\B-\overline{\B(1-\ve)}\right]$.

Note that for some $\de\in (0,\frac18 \ve)$ sufficiently small,
there exists  a compact embedded annulus of revolution
$E(\delta)\subset \left[(\overline{\B}-\B(1-\delta))\cap \B(p,\ve)\right]$ with
boundary circles in $\esf^2 \cup \left[\esf^2(1-\de)\right]$ and such  that the
radial projection $r(E(\de))\subset \esf^2$ is contained in the disk
$\B(p,\ve)\cap\esf^2$. Furthermore, $E(\delta)$ is also chosen to
have mean curvature greater than $H_0$ and with mean curvature
vector  outward pointing from the domain in
$\left[\overline{\B}-\B(1-\de)\right]-E(\de)$ which is contained in
$\B(p,\ve)$; for instance, $E(\de)$ can be chosen to be a piece of a
suitably scaled compact embedded annulus in some nodoid of constant
mean curvature one, see Figure~\ref{fig:3} Left. Assume now that 
$\ve$ is also chosen less than
$\frac{1}{10} d(x,y)$.

Assume that $n$ and $j$ are
chosen sufficiently large so that:
\ben
\item $f(A(n,j))\subset \left[\B - \overline{\B(1-\de)}\right]. $
\item Any circle in $\esf^2-\{x,y\}$ which represents the generator
of the first homology group $\HH_1 (\esf^2-\{x,y\})$ and whose
distance from $x$ and $y$ is at least $\de$, intersects the radial
projection $r(f(A(n,j)))\subset\esf^2$. This property holds since $x$
and $y$ are limit points of $f(A)$.
\item The radial projection of each of the two boundary curves of
$f(A(n,j))$ has diameter less than $\ve$. This condition is possible
to achieve since each of the components of $\partial f(A(n,j))$ has
image on either a disk component of $\cal W$, a face of $\cal F$ or
a component of $\partial \wt{N}(\G) - \cup \cal W$, and each of
these components and faces is  contained ambient balls of radius
$4\de$ by Lemma~\ref{size}, which in turn have radial projections of
diameter less than $\ve$. \een
\begin{figure}[htbp]
\begin{center} \includegraphics[width=.33\textwidth]{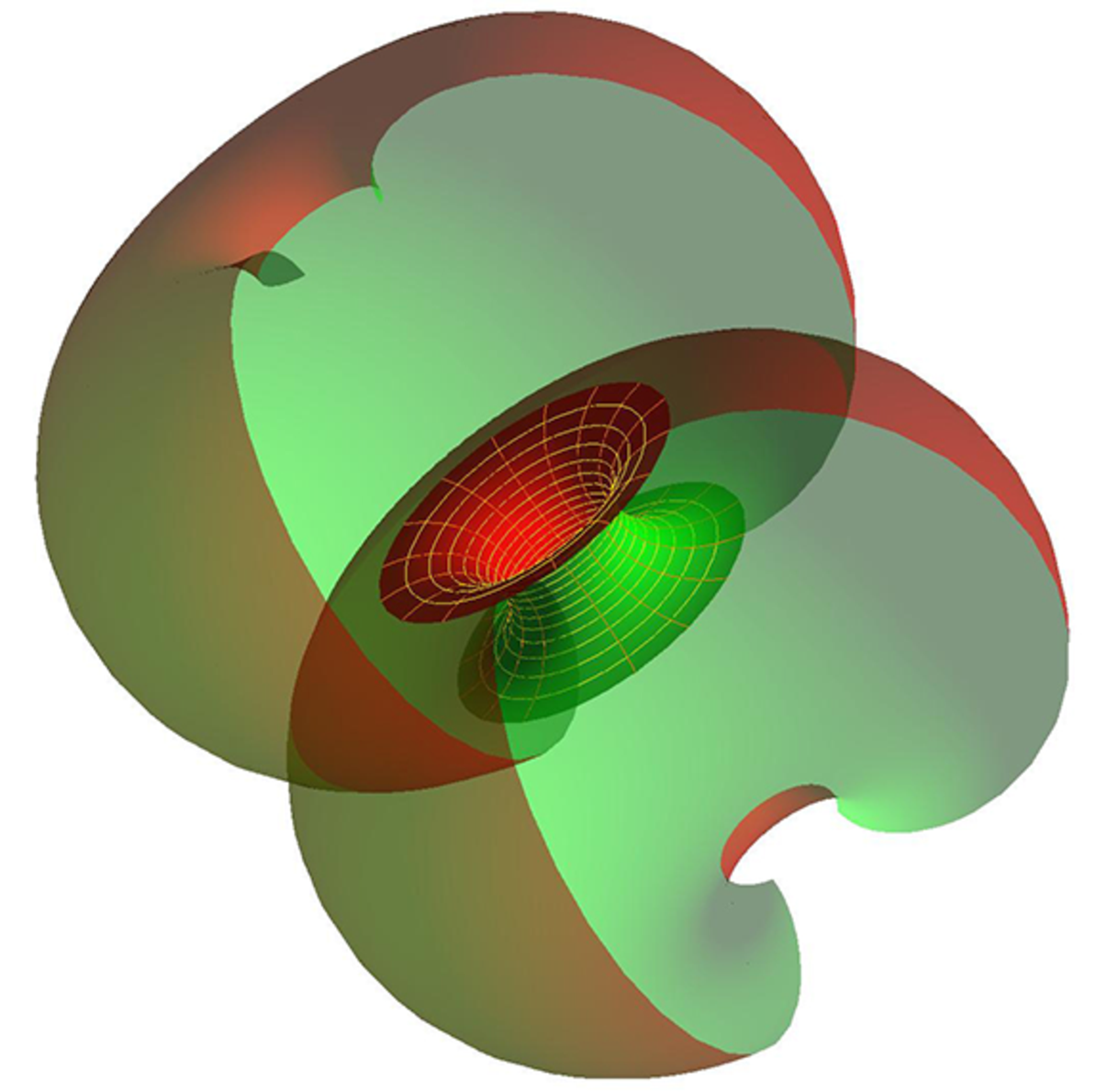}
\includegraphics[width=.6\textwidth]{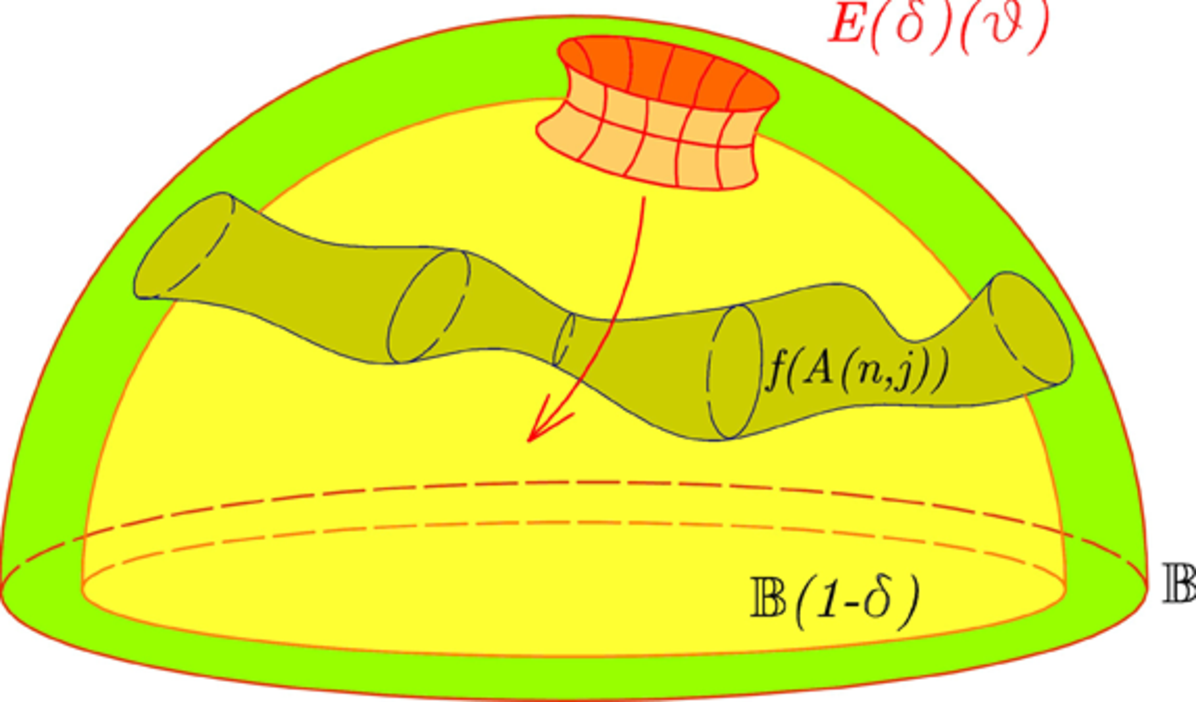}
\end{center}
\caption{Left: This figure shows a domain on the nodoid corresponding
to a scaling of $E(\de)$.  \, Right:  Since all of the $E(\de)(\vartheta)$ are disjoint
from $\partial f(A(n,j))$, a first point of contact in
$E(\delta)(\vartheta_0)\cap f(A(n,j))$ occurs an interior point
of $f(A(n,j))$.}\label{fig:3}
\end{figure}

By the above three properties and our choices of $\ve$ and $\de$, we
can choose a circle $S^1\subset \esf^2-\{x,y\}$ which intersects
$r(f(A(n,j)))$, and such that the $\ve$-neighborhood $S^1(\ve)$ of
$S^1$ is disjoint from the radial projection $r(\partial
f(A(n,j)))\subset \esf^2$ and each component of $\esf^2-r(S^1(\ve))$
contains points of $r(f(A(n,j)))$. Let $L$ be an oriented radial ray
which is an axis for the circle  $S^1$. For $\vartheta\in [0,2\pi)$,
consider the family of annuli $E(\delta)(\vartheta)$ obtained by
rotating  $E(\de)$ counterclockwise around $L$ by the angle
$\vartheta$. By elementary separation properties, there is a smallest
$\vartheta_0\in (0,2\pi)$ such that $E(\delta)(\vartheta_0)\cap f(A(n,j))\not
=\text{\O}$.  Since all of the $E(\de)(\vartheta)$ are disjoint from
$\partial f(A(n,j))$, a first point of contact in
$E(\delta)(\vartheta_0)\cap f(A(n,j))$ occurs an interior point of
$f(A(n,j))$, which must have absolute  mean curvature  on $A$ at
least equal to the minimum of the mean curvature of
$E(\delta)(\vartheta_0)$ (see Figure~\ref{fig:3} Right). But the mean curvature
of $E(\delta)(\vartheta_0)$ is greater than the absolute mean curvature
function of $A$. This contradiction completes the proof of
Lemma~\ref{l4.1}. \end{proof}

The next lemma follows immediately from the arguments presented in
the proof of Lemma~\ref{l4.1}; also see Figure~\ref{fig:3} Right.
We note that the constant $H_0$ in the statement of the next lemma is
the same constant which is the strict upper bound on the supremum of the
absolute mean curvature of $f\colon A\to \B$.
\begin{lemma} \label{lem:main} Given any $\ve\in(0,\frac{1}{4})$,
there exists an $\eta_0 \in (0,\frac{\ve}{10})$ that also depends on
$H_0$ such that the following statements hold. For any $\eta\in(0,\eta_0]$
and for any immersion $g\colon \Sigma \to \B-\B(1-\eta)$ of a compact surface
with boundary and absolute mean curvature less than $H_0$ such that
$g(\partial \Sigma)\subset \left[\overline{\B}(x,\eta)\cup \overline{\B}(y,\eta)\right]$
for two points $x,y\in \esf^2$ with $d(x,y)\geq \ve$, then either
$g(\Sigma)\subset \left[\overline{\B}(x,\eta)\cup \overline{\B}(y,\eta)\right]$
or $g(\Sigma)$ is $\ve$-close to every point in $\esf^2$. (Note that
it may be the case that $g(\partial \Sigma)$ is contained entirely
in one of the balls $\overline{\B}(x,\eta), \overline{\B}(y,\eta).$)
\end{lemma}

\begin{lemma}\label{lim:A(k,1)}
Given an $\ve \in (0,\frac{1}{2H_0})$, there exists an $n(\ve) \in \n$ such that for each
$k\geq n(\ve)$, there exists a point $y(k) \in \esf^2$ with
$f(A(k,1)) \subset \B(y(k), \ve)$ and $f(A(n(k)) \subset \left[\B -
\overline{\B}(1-\ve)\right]$.
\end{lemma}

\begin{proof} Fix $\ve \in (0,\frac{1}{2H_0})$ and let
$\eta=\min\{\eta_0,\eta_1\}$ where $\eta_0$ is given in Lemma~\ref{lem:main}
and depends on $\ve$ and $H_0$ and $\eta_1$ given in Lemma~\ref{lem:L2}.
Let $k(\eta)$ be the related integer given in Lemma~\ref{lem:chi}. We claim
that for $k\geq k(\eta)$, $f(A(k,1))\subset \B(z(k),\eta)$ and that $f(A(k(\eta))\subset\left[\B-\overline{\B}(1-\ve)\right]$, and so, by setting
$n(\ve)=k(\eta)$, this claim will complete the proof of the lemma.  By Lemma~\ref{lem:chi}, $f(A(k(\eta)))\subset\left[\B-\overline{\B}(1-\eta)\right]\subset \B-\overline{B}(1-\ve)$
and so it remains to
verify that $f(A(k,1))\subset \B(z(k),\eta)$.

Suppose that $f(A(k,1))$ contains a point outside of $\B(z(k),\eta)$. Let $K$ be
a nonempty component in $(f_{A(1,k)})^{-1}(\B-\B(z(k),\eta))$. By
Lemma~\ref{lem:maximum}, there is a point on $K$ which lies outside
of $\B(z(k), \frac{1}{H_0})$. Since $\ve\in (0,\frac{1}{2H_0})$ and
$\eta \leq \eta_0$, Lemma~\ref{lem:main} implies that the distance between  every point of
$\esf^2$ and $K$ is at most $\ve$. In particular, there exists a
point $p\in K$ such that $f(p)$ has distance greater than 1 from $z(k)$.

By the third statement in Lemma~\ref{lem:chi}, each boundary curve of
$K$ bounds a disk in $A(k,1)$. From the simple topology of an annulus
we find that exactly one boundary curve of $K$ bounds a disk
$D \subset A(k,1)$ and such that $K\subset D$. Next we apply
Lemma~\ref{lem:L2} to find an $X_i$-type curve $\beta_1\subset D$
which bounds a subdisk $D(\beta_1)$ which contains the point $p$
and which satisfies the other
properties in that lemma. In particular, we may assume the intrinsic
distance from $D(\beta_1)$ to $\partial D$ is at least $\eta$. By the
second statement in Lemma~\ref{lem:chi},
$f(\beta_1)\subset \B(y(\beta_1),\eta)$ for some point $y(\beta_1)\in \esf^2$.
By our previous arguments there exists a point $p_1\in D(\beta_1)$ such that
the distance from $f(p_1)$ to $y(\beta_1)$ is greater than one. So, we can apply
Lemma~\ref{lem:L2} again to obtain a subdisk $D(\beta_2)$,
$D\supset D(\beta_1)\supset D(\beta_2)$, where the intrinsic
distance from $\partial D(\beta_1)$ to $D(\beta_2)$ is at least $\eta$.
Repeating these arguments, induction gives the existence of a
sequence of disks $D\supset D(\beta_1)\supset\ldots\supset D(\beta_n)\supset\ldots$
such that the intrinsic distance from $D(\beta_n)$ to $\partial D$ is at least
$n \, \eta$. Since $D$ is compact, we obtain a contradiction which proves our
earlier claim that $f(A(k,1))\subset \B(z(k),\eta)$ for $k\geq k(\eta)$.
As we have already observed, this claim then proves the lemma.
\end{proof}

We now complete the proof of Theorem~\ref{th}. Fix some
$\ve'\in (0,\frac{1}{2H_0})$ and let $\eta_0\in (0,\frac{\ve'}{10})$
be the related number given in Lemma~\ref{lem:main}. Set $\ve=\eta_0$ and
let $n(\ve)$ be the integer given in Lemma~\ref{lim:A(k,1)}. By
Lemma~\ref{lim:A(k,1)}, for each
$i\in \N$,
$f(A(n(\ve),i))\subset \left[\B - \B(1-\eta_0)\right], f(\g_{n(\ve)})\subset \B(y(n(\ve)), \eta_0)$
and $f(\g_{n(\ve)+i})\subset \B(y(n(\ve)+i),\eta_0)$.
Since the limit set of $A(n(\ve))$ is all of $\esf^2$, there exists a
smallest $j\in \N$ such that  the distance between $y(n(\ve))$ and
$y(n(\ve)+j)$ is greater that $\ve'$. By Lemma~\ref{lem:main} and taking
into account that  $\B(y(n(\ve)),\eta_0)$ and $ \B(y(n(\ve)+j),\eta_0)$
are disjoint, then we conclude that
$f(A(n(\ve),j))$ must be $\ve'$ close to every point of $\esf^2$.

On the other hand, given $k \in \n$, $n(\ve) \leq k <n(\ve)+j$, we know
(by Lemma \ref{lim:A(k,1)}) that $f(A(k,1))\subset \B(y(k),\ve)$, for a 
suitable $y(k) \in\esf^2$. Moreover, the
choice of $j$ implies that $f(\gamma_k) \subset \B(y(n(\ve)),\ve'+\ve)$, for $k$ satisfying
$n(\ve) \leq k <n(\ve)+j$. So, by
 the triangle inequality we deduce $f(A(k,1)) \subset  \B(y(n(\ve)),\ve'+2 \ve)
 \subset \B(y(n(\ve)),2 \ve')$ for any $k$ satisfying $n(\ve) \leq k <n(\ve)+j$. 
 This implies $f(A(n(\ve),j))\subset \B(y(n(\ve)), 2\ve')$
which is impossible since $2\ve'<\frac{1}{10}$ and we have already seen that
$f(A(n(\ve),j))$ must be $\ve'$ close to every point of $\esf^2$. This contradiction
completes the proof of Theorem~\ref{th} in the case
$\overline{W}=\overline{\B}$.

For the general case where
$\overline{W}$ is a smooth  compact Riemannian manifold with
nonempty boundary,   small modifications of  the proof of
Theorem~\ref{th} in the special case $\overline{W}=\overline{\b}
\subset \rth$ also demonstrate that there exists a properly embedded
1-manifold $\Delta_W \subset W$, whose path components are smooth
simple closed curves, such that $\cd=W-\Delta_W$ is a Calabi-Yau
domain for any open surface with at least one annular end.
In carrying out these modifications in the smooth
compact 3-manifold $\overline{W}$,
it is convenient, to place the 1-manifold  $\Delta_{\overline{W}}$
in the union of small pairwise disjoint closed  $\ve$-neighborhoods
of the boundary components of $\overline{W}$ which have a natural
product structure derived from the distance function to the boundary
component.  The product structure simplifies the construction of the
related $1$-complex $\G_{\overline{W}}$ which has one component in
each of the $\ve$-neighborhoods of each boundary component of
$\overline{W}$. Also note that the
properness of any proper immersion of $A=\esf^1\times[0,\infty)$ into
$\cd$ guarantees that $A$ has a end representative which maps
into  the $\ve$-neighborhood of
exactly one of the boundary components of $\overline{W}$.
This discussion completes the proof of Theorem~\ref{th}.
\vspace{.2cm}

\begin{remark} {\rm The reader familiar with the paper~\cite{mmn}
might consider the question: Are the domains
$\D_{\cal F}\subset \rth$~\cite{mmn}, obtained
by removing a infinite  proper family $\F$ of horizontal circles
from $\B$, Calabi-Yau domains for surfaces with at least one annular end?
The answer to this question is no because for at least one such $\D_{\cal F}$
constructed in~\cite{mmn}, there exists a  proper, conformal,
complete embedding $f\colon \R^2\to \D$ with absolute mean curvature
function less than 1,  $f(\R^2)$ is a surface of revolution with
axis the $x_3$-axis, $f(\R^2)$ has intrinsic linear area growth and
has limit set $L(\R^2)=\esf^2$. The mean curvature function of
$f(\R^2)$ in this case contains points of mean curvature arbitrarily
close to 1 and also arbitrarily close to $-1$. In this case for
$\D$, the circles in $\F\subset \B$ are chosen to  have axis the
$x_3$-axis; the surface has the appearance of taking an infinite
connected sum of the spheres $\esf^2_k$, $k\in \n$, defined at the
beginning of Section~2,  joined by  small catenoidal type necks
centered along points along the $x_3$-axis which limit to the north
and south poles of $\esf^2$.}
\end{remark}

We conclude the paper with the following conjecture.
\begin{conjecture}
Let $\Delta \subset \B$ be the properly embedded one-manifold given in
the proof of Theorem~\ref{th}. If $\overline{B}$ is a smooth compact
Riemannian three-ball and $F\colon \overline{B} \to \B$ is a smooth
diffeomorphism, then $\D=B-F^{-1}(\Delta)$ is a Calabi-Yau domain
for {\bf any} noncompact  surface with compact boundary (possibly
empty). In particular, $\D=\B-\Delta$ does not admit any complete,
properly immersed open surfaces with bounded mean curvature.
\end{conjecture}

\center{Francisco Martín at fmartin@ugr.es\\
Department of Geometry and Topology, University of Granada, 18071 Granada, Spain}
\center{William H. Meeks, III at bill@math.umass.edu\\
Mathematics Department, University of Massachusetts, Amherst, MA
01003 USA}

\end{document}